\documentclass[preprint]{elsarticle}

\usepackage{lineno,hyperref}

\usepackage{amsmath}
\usepackage{amssymb}
\usepackage{mathtools}
\usepackage{tikz}
\usepackage{color}
\usepackage{graphicx}
\usepackage{bbold}
\usepackage{accents}

\modulolinenumbers[5]

\newcommand{\Ni}{\ensuremath{\mathcal{N}_i}}
\newcommand{\mk}{\ensuremath{<k>}}

\newcommand{\kmax}{\ensuremath{k_{\text{max}}}}
\newcommand{\tix}{\ensuremath{\theta_i^x}}
\newcommand{\tiy}{\ensuremath{\theta_i^y}}
\newcommand{\tjx}{\ensuremath{\theta_j^x}}
\newcommand{\tjy}{\ensuremath{\theta_j^y}}
\newcommand{\pix}{\ensuremath{\Psi_i^x}}
\newcommand{\piy}{\ensuremath{\Psi_i^y}}
\newcommand{\piz}{\ensuremath{\Psi_i^z}}
\newcommand{\dtix}{\ensuremath{\dot{\theta}_i^x}}
\newcommand{\dtiy}{\ensuremath{\dot{\theta}_i^y}}
\newcommand{\rpdt}{\ensuremath{r \psi_k d \theta}}
\newcommand{\spdt}{\ensuremath{\sigma \psi_k d \theta}}
\newcommand{\ga}{\ensuremath{\gamma}}
\newcommand{\al}{\ensuremath{\alpha}}

\newtheorem{thm}{Theorem}
\newtheorem{lem}[thm]{Lemma}
\newtheorem{cor}[thm]{Corollary} 
\newproof{pf}{Proof}

\journal{arXiv}

\begin{document}

\begin{frontmatter}

\title{Evolutionary Dynamics on a Regular Networked Structured and Unstructured Multi-population
}

\author[rug]{Wouter Baar\corref{mycorrespondingauthor}}
\cortext[mycorrespondingauthor]{Corresponding author}
\ead{w.baar@rug.nl}

\author[rug,ita]{Dario Bauso}
\ead{d.bauso@rug.nl}

\address[rug]{Jan C. Willems Center for Systems and Control, ENgineering and TEchnology institute Groningen (ENTEG), Faculty of Science and Engineering, University of Groningen, 9747 AG Groningen, The Netherlands}
\address[ita]{Dipartimento di Ingegneria, Universit\`a di Palermo, 90128 Palermo}

\begin{abstract}
In this paper we study collective decision making on a multi-population, represented by a regular network of groups of individuals. Each group consists of a collection of players and every player can choose between two options. A group is characterised by variables denoting the fractions of individuals committed to each respective option, and they are influenced by the state of neighboring groups. First, we study its steady-state and show that the equilibrium is a consensus equilibrium. We also derive a sufficient condition for local asymptotic stability. Then, we study a structured model where every population is now assumed to represent a structured complex network. We conclude the paper with simulations, corroborating the obtained theoretical findings.
\end{abstract}


\end{frontmatter}

\linenumbers

\section{Introduction}
In this paper we study networked bio-inspired models, where a group of individuals has the objective of reaching consensus on one of the available options. In recent years there has already been a surge in amount of published literature on this topic, see for example \cite{dario}, \cite{j13}, \cite{hill} and \cite{bauso3}. Usually, a group of individuals is considered and every player of the group can choose between two options, $1$ or $2$, or can choose to be uncommitted. Rather than studying individual behavior of each player, we are only interested in how many players choose a certain strategy. To do this, we model the mean-field approximation, by taking the population size very large. 

In the existing literature, often the dynamics of a single group of individuals is studied, \cite{dario}, \cite{j13}, \cite{bauso3}. Here, the evolutionary dynamics are inspired on a biological model on swarms of bees, see e.g. \cite{bee1}, \cite{bee2}, \cite{bee3}. Networked bio-inspired evolutionary models also arise in consensus networks \cite{multipoppaper},  \cite{opiniondyn}, \cite{hill}, \cite{ming}. Here, we can think of groups of individuals that need to decide between a left or right political party, and the cross-inhibitory signals can then be thought of as smear campaigns by the opposing team. Finally, networked systems can be used to study diffusion models. For example, in \cite{bullogoed}, they are used to study how a disease spreads over a population. Individuals are now either susceptible or infected. If there are two competing viruses spreading over the population, we talk about a bi-virus model \cite{basar}.

A preliminary version of this paper has appeared as a conference proceeding \cite{onzeconf}. The additional main contributions are the study of the case with structured environment, the addition of a lemma proving the well-posedness of the problem, and extended proofs of the obtained results.

As an element of novelty in comparison with the existing literature, in this paper we will deal with evolutionary dynamics on a multi-population, as opposed to a single population. This networked multi-population is represented by a graph, where each node represents a population. Every population is characterized by three variables, the fractions of players that are committed to option 1, committed to option 2, and uncommitted. These fractions depend on the state of the neighboring populations. We assume throughout the paper that the graph of the networked multi-population is unweighted and regular, by which we mean that every node has the same amount of neighbors.

\section{The Unstructured Networked Model}\label{sec:model}
Consider a multi-population that consists of multiple groups of players, where it is assumed that players within a group can all interact with each other. The evolutionary dynamics are now as follows. Every player can spontaneously decide to commit to an option (parametrized by $\gamma$), or can decide to commit to an option by means of imitation (parametrized by $r$), because we assume the individuals are crowd-seeking, and thus are attracted to the option that has the most committed players. Players can also decide to leave their choice, by either spontaneous abandonment (parametrized by $\alpha$), or by being lured to become uncommitted by other players. We then also talk about cross-inhibitory signals that are sent to opposing players to attract them to become uncommitted (parametrized by $\sigma$). Finally, when modeling the dynamics of such systems, the states of the model are $x,y,z$, that denote the fraction of committed individuals to option 1, the fraction of committed individuals to option~$2$ and the fraction of uncommitted individuals, respectively. We introduce a subscript $i$ to denote the $i^{\text{th}}$ population. Note that we have $x_i + y_i + z_i = 1$ at all times, or equivalently, $z_i =1-x_i -y_i$. This enables us to describe the dynamics solely in terms of $x_i$ and $y_i$: 
\begin{equation}\small 
\begin{split} \label{eq:regulmicrodyn}
  \dot{x}_i &= \Big( \gamma + r \sum\limits_{j \in \Ni}  x_j \Big) (1 - x_i - y_i )  - x_i \Big( \alpha + \sigma \sum_{j \in \Ni}  y_j \Big) , \\ 
  \dot{y}_i &= \Big( \gamma + r \sum_{j \in \Ni}  y_j \Big)  (1 - x_i - y_i )  - y_i \Big( \alpha + \sigma \sum_{j \in \Ni}  x_j \Big) ,  
 \end{split} 
\end{equation} \normalsize 
\begin{figure}[t!]
    \centering
    \begin{tikzpicture}[scale=0.8]
    \node (1) at (0,0) [shape= circle, draw=black, label=below:{$z_i$}] {}; 
    \node (2) at (-4, 0) [shape=circle, draw = black, label=below:{$x_i$}] {}; 
    \node (3) at (4,0) [shape=circle, draw=black, label=below:{$y_i$}] {}; 
    \draw (1) edge[->, thick, bend right=10] node[above] {$\gamma + r \sum  x_i$} (2); 
    \draw (2) edge[->, thick, bend right=10] node[below] {$\alpha + \sigma \sum y_i$} (1); 
    \draw (1) edge[->, thick, bend left=10] node[above] {$\gamma + r \sum y_i$} (3); 
    \draw (3) edge[->, thick, bend left=10] node[below] {$\alpha + \sigma \sum x_i$} (1); 
    \end{tikzpicture}
    \caption{The evolutionary dynamics at the node level.}
    \label{fig:ourmodel}
\end{figure}
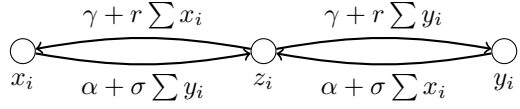
\noindent and of course $\dot{z}_i = - \dot{x}_i - \dot{y}_i$. These dynamics and the role of the parameters are summarized in Figure \ref{fig:ourmodel}. Throughout it is assumed that we consider a regular and unweighted graph, so the cardinality of the neighbor set $| \Ni |$ is the same for every node. The set of feasible states $\Delta$, defined by 
\begin{align*}
    \Delta = \{ (x,y) \in \mathbb{R}^{2n} \mid x_i \geq 0, y_i \geq 0, x_i + y_i \leq 1 \} ,  
\end{align*} 
is positively invariant. This makes the problem we study in this paper well-posed. 
\begin{lem}
For any initial condition $(x(0), y(0)) \in \Delta$, the state will remain in this set, $(x(t), y(t)) \in \Delta$, for all $t>0$.
\end{lem}
\begin{pf}
This is a corollary from Nagumo's Theorem \cite{frabl}. To see this, we compute the direction of the vector field of the points on the boundary of $\Delta$, and show that they point inwards the set. Consider any node $i$. Recall that $\gamma, r , \alpha, \sigma$ are all positive valued parameters, and assume that $x_j, y_j \geq 0$ for all $j \neq i$. We now show that $\Delta_i = \{ (x_i, y_i ) \mid x_i \geq 0, y_i \geq 0, x_i + y_i \leq 1 \}$ is a positively invariant set. 
\mbox{The boundary of $\Delta_i$ is}
\begin{align*}
    \partial \Delta_i &= \{ (x_i, y_i ) \in \Delta_i \mid y_i = 0 \} \cup \{ (x_i, y_i ) \in \Delta_i \mid x_i = 0 \} \\ &\quad \cup \{ (x_i, y_i ) \in \Delta_i \mid  x_i + y_i = 1 \} \\ &=: \partial \Delta_i^1 \cup \partial \Delta_i^2 \cup  \partial \Delta_i^3 . 
\end{align*}
For the points in $\partial \Delta_i^1$, we have that 
\begin{align*}
  \dot{y}_i &= \Big( \gamma + r \sum\nolimits_{j \in \Ni} a_{ij} y_j \Big)  (1 - x_i ) \geq 0  ,
\end{align*}
so these points are moving upwards inside $\Delta_i$. For points in $\partial \Delta_i^2$ we have that they point inside $\Delta_i$ as well, since 
\begin{align*}
     \dot{x}_i &= \Big( \gamma + r \sum\nolimits_{j \in \Ni} a_{ij} x_j \Big) ( 1  - y_i ) \geq 0 .  
\end{align*}
And finally, for points in $\partial \Delta_i^3$ we have that
\begin{align*}
    \dot{x}_i &= - x_i \Big( \alpha + \sigma \sum\nolimits_{j \in \Ni} a_{ij} y_j \Big) \leq 0 , \\
    \dot{y}_i &=  - y_i \Big( \alpha + \sigma \sum\nolimits_{j \in \Ni} a_{ij} x_j \Big) \leq 0 , 
\end{align*}
which points inward $\Delta_i$ as well. To conclude, we see that for any point on the boundary of $\Delta$, the vector field points inside $\Delta$. As a consequence of Nagumo's Theorem, we can conclude that the set $\Delta$ is positively invariant. 
\end{pf}
We are now ready to establish the following result.  
\begin{thm}\label{thm:conseq} 
The system described by \eqref{eq:regulmicrodyn} admits a consensus equilibrium $(x^*, y^*, z^* ) = (\xi  \mathbb{1}_n, \mu \mathbb{1}_n , \zeta \mathbb{1}_n )$  with 
\begin{itemize}
    \item Case 1. If $\xi = \mu$, 
    \begin{equation*}
        \begin{split}
            \xi = \mu = \tfrac{- (2 \gamma - rd + \alpha) + \sqrt{(2 \gamma - rd + \alpha )^2 + 4 \gamma (2 rd + \sigma d)  }}{2 ( 2rd + \sigma d)} ,\\
            \zeta = 1 - \tfrac{- (2 \gamma - rd + \alpha) + \sqrt{(2 \gamma - rd + \alpha )^2 + 4 \gamma (2 rd + \sigma d)  }}{ ( 2rd + \sigma d)} .
        \end{split}
    \end{equation*}
    \item Case 2. If $\zeta = \frac{\alpha }{dr}$, 
\begin{equation*}
\begin{split}
    \xi = \tfrac{- (\frac{ \alpha}{r} -  d) \pm \sqrt{(\frac{\alpha}{r}- d)^2 -  \frac{4  \gamma \alpha}{r \sigma }}}{2  d} ,\ \   \mu = 1-  \xi  - \tfrac{\alpha}{r d} , \ \ \zeta = \tfrac{\alpha}{ rd } .
    \end{split}
\end{equation*} 
\end{itemize}
\end{thm}
Our notion of consensus refers to the fact that the fractions of (un)commitment, $\xi, \mu,$ and $\zeta$, are the same for every population.

\begin{pf}
Let the consensus equilibrium be given by $(x^*, y^*, z^* ) = (\xi  \mathbb{1}_n, \mu \mathbb{1}_n , \zeta \mathbb{1}_n )$, with $\zeta = 1 - \xi - \mu $. At this equilibrium, and since $\dot{x}_i=0$ and $\dot{y}_i=0$ at an equilibrium, the equations in \eqref{eq:regulmicrodyn} reduce to 
\begin{equation}\begin{split} \label{eq:regreduc} 
      0 = \dot{x_i} = ( \gamma + r d \xi ) \zeta - \xi ( \alpha + \sigma d \mu ) ,\\
    0 = \dot{y_i} = ( \gamma + r d \mu ) \zeta - \mu ( \alpha + \sigma d \xi ) . \end{split}
\end{equation} 
From \eqref{eq:regreduc} it also follows that
\begin{align*}
    0 = \dot{x}_i - \dot{y}_i = rd ( \xi - \mu ) \zeta - \xi \alpha + \mu \alpha
    = (rd \zeta - \alpha ) ( \xi - \mu ) . 
\end{align*}
From this expression we obtain that either $\xi = \mu$ or $\zeta = \tfrac{\alpha}{rd}$, so we distinguish two cases. Note that these two cases are not mutually exclusive. In the first case, if $\xi = \mu$, the first equation of \eqref{eq:regreduc} reduces to
    \begin{align*}
    0 = \dot{x}_i  &= ( \gamma + r d \xi ) ( 1 - 2 \xi ) - \xi ( \alpha + \sigma d \xi )\\ &=  
    \gamma - 2 \gamma \xi + r d \xi - 2 r d \xi^2 - \alpha \xi - \sigma d \xi^2 \\  &=
    (2 r d + \sigma d ) \xi^2 + ( 2 \gamma - rd + \alpha) \xi - \gamma .  
\end{align*} 
This is a quadratic equation in $\xi$ and the solution is
\begin{align*}
    \xi = \tfrac{- (2 \gamma - rd + \alpha) \pm \sqrt{(2 \gamma - rd + \alpha )^2 + 4 \gamma (2 rd + \sigma d)  }}{2 ( 2rd + \sigma d)} . 
\end{align*}
Note that in general one obtains two solution to a quadratic equation, however the solution with a minus sign in front of the square root cannot be an equilibrium since in that case the value of $\xi$ would be strictly negative. 
Finally, $\mu = \xi$ and $\zeta = 1 - 2 \xi$. Secondly, if $\zeta = \tfrac{\alpha}{rd}$, using the fact that $\mu = 1 - \xi - \tfrac{\alpha}{rd}$ now, we obtain 
\begin{align*}
    0 = \dot{x}_i &= ( \gamma + r d \xi) \tfrac{\alpha}{r d} - \xi ( \alpha + \sigma d (1 - \xi - \tfrac{\alpha}{rd}) )  \\ &=
    d \xi^2 + ( \tfrac{\alpha}{r} - d) \xi + \tfrac{\gamma \alpha }{r d \sigma }.
\end{align*}
This is a quadratic equation in $\xi$ and the solution is
\[   \xi = \tfrac{- (\frac{ \alpha}{r} -  d) \pm \sqrt{(\frac{ \alpha}{r}- d)^2 -  \frac{4  \gamma \alpha}{r \sigma }}}{2  d} . \] 
It can happen that for some values of parameters both solutions are within the interval $[0,1]$. In that case we do not have a unique equilibrium. This completes the proof. 
\end{pf}
For the consensus equilibrium established in Theorem \ref{thm:conseq}, we have the following sufficient condition for exponential stability. We will show that two inequalities must hold on the strength of the cross-inhibitory signals, and these conditions can be checked a priori. 
\begin{thm}
The consensus equilibrium $x^* = (\xi \mathbb{1}_n, \mu \mathbb{1}_n, \zeta \mathbb{1}_n)$ that is reached under dynamics \eqref{eq:regulmicrodyn} is locally exponentially stable if the following inequalities on the cross-inhibitory signal $\sigma$ hold: 
\begin{align*}
    \sigma > \dfrac{\alpha - r d (1 - \xi - \mu)}{d (1 - \mu)} \ \text{ and } \   \sigma > \dfrac{\alpha - rd (1 - \xi - \mu )}{d ( 1 - \xi)}. 
\end{align*}
\end{thm}
\begin{pf} 
We compute the Jacobian and evaluate it at the equilibrium $x^* = (\xi \mathbb{1}_n, \mu \mathbb{1}_n, \zeta \mathbb{1}_n )$. Consider the dynamics described by \eqref{eq:regulmicrodyn}, then the partial derivatives of $\dot{x}_i$ are given by the following expressions
\begin{align*}
\begin{array}{ll}
    \begin{array}{l} \dfrac{\partial \dot{x}_i }{\partial x_i} = - \gamma - r \sum_{j \in \Ni} x_j \\ \qquad \quad -   \alpha - \sigma \sum_{j \in \Ni} y_j, \end{array} &  \dfrac{\partial \dot{x}_i}{\partial y_i} = - \gamma - r \sum_{j \in \Ni} x_j , \\  
    \dfrac{\partial \dot{x}_i }{\partial x_j} = \begin{cases} 0 \text{ if } j \notin \Ni, \\ r z_i  \text{ if } j \in \Ni,\end{cases}  &     \dfrac{\partial \dot{x}_i}{\partial y_j} = \begin{cases} 0 \text{ if } j \notin \Ni, \\ - \sigma x_i \text{ if } j \in \Ni. \end{cases}  \end{array} 
\end{align*}
The partial derivatives of $\dot{y}_i$ are obtained in a similar fashion. 
We evaluate the expressions of the partial derivatives at the equilibrium $x^* = (\xi \mathbb{1}_n, \mu \mathbb{1}_n, \zeta \mathbb{1}_n)$ and we make use of the fact that $1 - \xi - \mu = \zeta$. This gives the following Jacobian 
\begin{align}\label{eq:jacobian}
    J & = \left. \begin{bmatrix} \dfrac{\partial (\dot{x}_1, \cdots, \dot{x}_n , \dot{y}_1, \cdots, \dot{y}_n )}{ \partial (x_1, \cdots, x_n, y_1, \cdots, y_n)}  \end{bmatrix} \right|_{x=x^*} \nonumber  \\ 
    &= \left. \begin{bmatrix} \begin{array}{c|c} \dfrac{\partial ( \dot{x}_1 , \cdots , \dot{x}_n ) }{\partial ( x_1 , \cdots , x_n ) } & \dfrac{\partial (  \dot{x}_1 , \cdots , \dot{x}_n ) }{\partial ( y_1 , \cdots , y_n )} \\ \hline  \dfrac{\partial ( \dot{y}_1 , \cdots , \dot{y}_n ) }{\partial ( x_1 , \cdots, x_n )} & \dfrac{\partial ( \dot{y}_1 , \cdots , \dot{y}_n ) }{\partial ( y_1, \cdots , y_n ) } \end{array}  \end{bmatrix} \right|_{x=x^*} \nonumber \\ 
    &= \begin{bmatrix}\begin{array}{c|c} ( r - r\xi - r \mu ) A - & - \sigma A - (\gamma + rd \xi) I_n \\ ( \gamma + rd \xi + \alpha + \sigma d \mu ) I_n  & \\ \hline   & (r - r \xi - r \mu) A -  \\ - \sigma A - (\gamma + r d \mu ) I_n & (\gamma + r d \mu + \alpha + \sigma d \xi ) I_n  \end{array} \end{bmatrix} \nonumber , 
\end{align}
where $A$ denotes the adjacency matrix of the network. To prove exponential stability the eigenvalues of the above matrix  should be contained in the open left half of the complex plane. Computing the eigenvalues of the above $2n \times 2n$ Jacobian matrix is rather difficult, however, we can find an estimate using the Gershgorin circle theorem. Gershgorin circle theorem tells us that each eigenvalue $\lambda_i$ is contained in a circle around $J_{ii}$ with radius $\sum_{i \neq j} | J_{ij} |$. We note that all diagonal entries $J_{ii}$ are strictly negative, so if we require that the radius of the circle does not cross the imaginary axis we know that the eigenvalue has a strictly negative real part. Using this, we find a sufficient condition for stability. So for $i= 1,\ldots, n$ we require
\begin{align*} 
- J_{ii} &> \sum_{j \neq i } | J_{ij} | \\ 
    \gamma + rd \xi + \alpha + \sigma d \mu &>  r d \zeta + \sigma d + \gamma + r d \xi \\ 
     \alpha + \sigma d \mu &> r d \zeta + \sigma d \\
    \sigma &> \dfrac{r d \zeta - \alpha }{d  \mu - d} = \dfrac{ \alpha - r d \zeta }{d (1 - \mu) }. 
\end{align*}
Analogously, for $i = n+1, \ldots, 2n$ we must have $- J_{ii} > \sum_{j \neq i } | J_{ij} |$ and similar computations show that this gives the second constraint.
This completes the proof.
\end{pf}

\section{A Structured Networked Multi-population}\label{sec:struc} 
So far, we did not impose any structure on the individual populations. A single population was merely characterised by three variables, namely the fraction of committed individuals to option 1, the fraction of committed individuals to option~2 and the fraction of uncommitted individuals. The internal structure of the population was not considered. That will change in the remainder of this paper. From now on, we will assume that every group of the multi-population is characterized by a structured environment \cite{thepaperonmydesk}, \cite{bauso3}. The key idea here, is that individuals from a population with a high connectivity to other players (in the same population) can have different inclinations to commit to option 1 or 2 than compared to a player with a low connectivity. 

To make this more precise, we will now study the case where each population of the networked multi-population is represented by a complex network. A complex network is a network on a large number of nodes, where the node degree distribution follows a power law. We let $P(k)$ be the distribution of the node degrees, for $k=1,2,\ldots,k_{\text{max}}$. Then we denote $x_i^k$ as the portion of the $i^{\text{th}}$ population with $k$ connections, that are committed to the first option. Similarly for $y_i^k$ and $z_i^k$, and we have $z_i^k = 1 - x_i^k - y_i^k$. We also refer to the pair $(x_i^k, y_i^k)$ (or the triple $(x_i^k, y_i^k, z_i^k)$) as the $k^{\text{th}}$ cluster of the $i^{\text{th}}$ population. It is noteworthy that the full state of the networked complex network population has as state vectors 
\begin{align*}
    \left( x_1^1 , \  x_1^2 , \  \hdots , \ x_1^{k_{\text{max}}} , \  \hdots , \  x_n^1 , \  x_n^2 , \  \hdots , \  x_n^{k_{\text{max}}} \right)^T , \\  \left( y_1^1 , \  y_1^2 , \ \hdots , \  y_1^{k_{\text{max}}} , \  \hdots , \  y_n^1 , \  y_n^2 , \  \hdots , \  y_n^{k_{\text{max}}} \right)^T , 
\end{align*}
where we assumed that for every complex network, the degree distribution $P(k)$ is the same, and hence the maximal node degree $k_{\text{max}}$ is also equal for every group. Let us now consider the dynamics for a single cluster from a single complex network of the network. In other words, we will now come up with dynamics for $x_i^k$ and $y_i^k$. We introduce $\psi_k = \tfrac{k}{k_{\text{max}}}$ as a parameter that captures the connectivity of a cluster. We have that $\psi_k$ is close to zero if $k$ is very small (no connectivity) and $\psi_k$ is close to one if $k$ is close to $k_{\text{max}}$ (full connectivity). For the discrete random variable $k$, with probabilities $P(k)$, we denote the mean by $<k>$ with 
\[ <k> \  = \sum\nolimits_{k=1}^{k_{\text{max}}} k \cdot P(k)  . \]
Next up, we introduce $\theta_i^x$ and $\theta_i^x$ as the probabilities that a link randomly chosen will point to a player that uses strategy $x$ or $y$ in group $i$, respectively \cite{thetaref}. These new variables capture the first moment. They are defined as 
\begin{align*}
     \tix = \dfrac{1}{\mk} \sum_{k = 1}^{k_{\text{max}}} k \cdot P(k) \cdot x_i^k , \ \ 
    \tiy  = \dfrac{1}{\mk} \sum_{k = 1}^{k_{\text{max}}} k \cdot P(k) \cdot y_i^k . 
\end{align*}
The dynamics of cluster $k$ in group $i$ are now given by \small
\begin{align}\begin{split} \label{eq:micromod} 
    \dot{x}_i^k &= ( 1 - x_i^k - y_i^k ) \Big( r \psi_k \sum_{j \in \mathcal{N}_i} \theta_j^x + \gamma  \Big) - x_i^k \Big( \alpha + \sigma \psi_k \sum_{j \in \mathcal{N}_i} \theta_j^y \Big) , \\ 
    \dot{y}_i^k &= ( 1 - x_i^k - y_i^k ) \Big( r \psi_k \sum_{j \in \mathcal{N}_i} \theta_j^y + \gamma  \Big) - y_i^k \Big( \alpha + \sigma \psi_k \sum_{j \in \mathcal{N}_i} \theta_j^x \Big) . \end{split} 
\end{align} \normalsize 
Note that these dynamics are normalised, in the sense that $x_i^k + y_i^k + z_i^k = 1$ for every $k \in \{1,2,\ldots, k_{\text{max}} \}$ and for every $i \in \{ 1,2, \ldots, n \}$. We made this normalisation because in this case, $z_i^k = 1 - x_i^k - y_i^k$, and this variable became redundant. These dynamics are similar to the model \eqref{eq:regulmicrodyn} except we now sum over the multiple $\theta$ parameters, and the connectivity parameter $\psi_k$ now play a role as well, since we are looking at a multi-population of complex networks. By aggregating over the differential equations of $\dot{x}_i^k$ and $\dot{y}_i^k$ using $\tfrac{1}{<k>} \sum_k k P(k)$, we obtain dynamics for \dtix and \dtiy, which are given by 
\begin{align*}
    \dtix &= (- \gamma - \alpha ) \tix - \gamma \tiy + \tfrac{r \Psi_i^z}{ \kmax} \sum_{j \in \Ni} \tjx - \tfrac{\sigma \pix}{\kmax} \sum_{j \in \Ni} \tjy + \gamma , \\
    \dtiy &= - \gamma \tix + (- \gamma - \alpha) \tiy -   \tfrac{ \sigma \piy}{\kmax} \sum_{j \in \Ni} \tjx + \tfrac{r \Psi_i^z }{\kmax} \sum_{j \in \Ni} \tjy + \gamma , 
\end{align*}
where we have simplified our notation by introducing 
\begin{align*}
\begin{array}{ll} 
    V = \sum\nolimits_k k^2 P(k), & \pix = \tfrac{1}{\mk} \sum\nolimits_k k^2 P(k) x_i^k , \\
    \piy = \tfrac{1}{\mk} \sum\nolimits_k k^2 P(k) y_i^k , &  \Psi_i^z = \tfrac{1}{\mk} \sum\nolimits_k k^2 P(k) z_i^k . \end{array} 
\end{align*}
Note that we have $V = = \mk \left( \pix + \piy + \Psi_i^z \right)$. These variables capture the second moment. We can write the dynamics for $\tix$ and $\tiy$ for $i=1,2,\ldots,n$ in vector form, by introducing 
\begin{align*}
    \theta^x = \left( \theta_1^x, \  \theta_2^x, \  \ldots , \  \theta_n^x \right)  , \ \ \theta^y = \left(  \theta_1^y , \theta_2^y , \ \ldots , \ \theta_n^y \right)  . 
\end{align*} 
We then have  \small 
\begin{align} 
\begin{split} \label{eq:micmac}
    \begin{pmatrix} \dot{\theta}^x \\ \dot{\theta}^y \end{pmatrix} &= \begin{bmatrix} (- \gamma - \alpha ) I + \tfrac{r}{\kmax} \Psi^z A & - \gamma I - \tfrac{\sigma}{\kmax} \Psi^x A \\ - \gamma I - \tfrac{\sigma}{\kmax} \Psi^y A & (-\gamma - \alpha) I + \tfrac{r}{\kmax} \Psi^z A \end{bmatrix} \begin{pmatrix} \theta^x \\ \theta^y \end{pmatrix} \\ &\quad + \begin{pmatrix} \gamma \mathbb{1} \\ \gamma \mathbb{1} \end{pmatrix}  , 
    \end{split} 
\end{align} \normalsize 
where we use $A$ to denote the adjacency matrix of the multi-population, and where we use $\Psi^x. \Psi^y$ and $\Psi^z$ to denote  
\begin{align*}
    \Psi^x = \text{diag}(\Psi_1^x , \Psi_2^x , \ldots, \Psi_n^x) ,\\
    \Psi^y = \text{diag}(\Psi_1^y , \Psi_2^y , \ldots, \Psi_n^y)  ,\\
    \Psi^z = \text{diag}(\Psi_1^z , \Psi_2^z , \ldots, \Psi_n^z) .
\end{align*}
Equilibria of the above  model are in general hard to find, but the equilibrium $(\theta^{x*}, \theta^{y*})$ satisfies at least the expression \small 
\begin{align*} 
    \begin{pmatrix} \theta^{x*} \\ \theta^{y*} \end{pmatrix} &=   \begin{bmatrix} ( \gamma + \alpha ) I - \tfrac{r}{\kmax} \Psi^z A &  \gamma I + \tfrac{\sigma}{\kmax} \Psi^x A \\  \gamma I + \tfrac{\sigma}{\kmax} \Psi^y A & (\gamma + \alpha) I - \tfrac{r}{\kmax} \Psi^z A \end{bmatrix}^{-1} \begin{pmatrix} \gamma \mathbb{1} \\ \gamma \mathbb{1} \end{pmatrix}  , 
\end{align*} \normalsize 
where the matrix is nonsingular if it is diagonally dominant. Computation of an equilibrium analytically is a difficult task, and is beyond the scope of this paper. Suppose for now, an equilibrium $(\theta^{x*}, \theta^{y*})$ is found. Then we can conclude the following result on stability. 
\begin{thm}\label{thm:cor9}
Let the multi-population be characterized by a regular and unweighted graph with degree $d$ and with adjacency matrix $A = (a_{ij})$. Then, an equilibrium $(\theta^{x*}, \theta^{y*})$ of \eqref{eq:micmac} is locally asymptotically stable if the following pair of inequalities are satisfied. 
\begin{align*}
    \alpha &> \tfrac{r}{\kmax} \Psi_*^z d + \tfrac{\sigma}{\kmax} \Psi_*^x d , \\
    \alpha &> \tfrac{\sigma }{\kmax} \Psi_*^y d  + \tfrac{r}{\kmax} \Psi_*^z d , 
\end{align*} 
where $\Psi_*^x, \Psi_*^y, \Psi_*^z$ denote the second moments $\pix, \piy, \piz$, since the latter are the same for every population $i$.
\end{thm} 
It should be noted that the equilibrium values play a role in the computation of $\Psi_*^x, \Psi_*^y$ and $\Psi_*^z$.   
\begin{pf}
The result follows directly from applying Greshgorin's circle theorem to the affine system of \eqref{eq:micmac}. Computing the eigenvalues of the system matrix explicitly is quite difficult, however they can be estimated by using Gershgorin's circle theorem. For every row of the system, since the diagonal entry is strictly negative, requiring that it is larger in magnitude than the sum of the off-diagonal entries will tell us that the eigenvalues are contained in the open left half complex plane. So for $i=1,2,\ldots,n$ we have that 
\begin{align*}
    |  - \gamma - \alpha | &>  \left|  \tfrac{r}{\kmax} \piz  \sum_{j = 1 }^n  a_{ij} \right| + | - \gamma |  + \left| -  \tfrac{\sigma}{\kmax} \pix \sum_{j =1}^n a_{ij } \right|  \\
    \alpha &> \tfrac{r}{\kmax} \piz \sum_{j=1}^n a_{ij} + \tfrac{\sigma}{\kmax} \pix \sum_{j=1}^n a_{ij} \\ 
    \alpha &> \tfrac{r}{\kmax} \piz d + \tfrac{\sigma}{\kmax} \pix d,  
\end{align*}
since $\sum_{j=1}^n a_{ij} = | \Ni | = d$ for all $i$. Similarly, for rows $i=n+1,\ldots,2n$ we require that 
\begin{align*}
    \alpha &> \tfrac{r }{\kmax} \piz d + \tfrac{\sigma}{\kmax} \piy d . 
\end{align*}
Furthermore we note that since we consider a regular graph, we reach a consensus equilibrium by which we mean that the values for $x_i^k$ and $y_i^k$ are the same for every $1 \leq i \leq n$, as in Theorem \ref{thm:conseq}. We thus have 
\begin{align*}
    \pix &= \tfrac{1}{\mk} \sum_k k^2 P(k) x_i^k = \tfrac{1}{\mk} \sum_k k^2 P(k) \xi^k = \Psi_*^x  , \\
    \piy &= \tfrac{1}{\mk} \sum_k k^2 P(k) y_i^k = \tfrac{1}{\mk} \sum_k k^2 P(k) \mu^k = \Psi_*^y , \\
    \piz &= \tfrac{1}{\mk} \sum_k k^2 P(k) z_i^k = \tfrac{1}{\mk} \sum_k k^2 P(K) \zeta^k = \Psi_*^z .  
\end{align*}
This completes the proof. 
\end{pf} 
In the special case of reaching a symmetric equilibrium as in case 1 of Theorem \ref{thm:conseq}, we have the following corollary. 
\begin{cor} 
Let the underlying network of the multi-population be a regular and unweighted graph, with degree $d$. Then, any symmetric equilibrium of \eqref{eq:micmac} is locally asymptotically stable if the following condition on the cross-inhibitory signal holds
\begin{align}\label{eq:stabcond} 
    \sigma < 2 r - \dfrac{r V}{\mk \Psi_*} + \dfrac{\alpha \kmax}{\Psi_* d} , 
\end{align}
where $\Psi_* = \Psi_*^x = \Psi_*^y$. 
\end{cor}
\begin{pf}
Since we have a symmetric equilibrium, $\Psi_*^x = \Psi_*^y = \Psi_*$. Furthermore this means $\Psi_*^z = \tfrac{V}{\mk} - \Psi_*^x - \Psi_*^y = \tfrac{V}{\mk} - 2 \Psi_*$. Starting from Theorem \ref{thm:cor9} and performing some algebraic manipulations, we have 
\begin{align*} 
     \alpha &> \tfrac{r}{\kmax} \Psi_*^z d + \tfrac{\sigma}{\kmax} \Psi_*^x d \\
    \alpha &> \tfrac{r}{\kmax} \left( \tfrac{V}{\mk} - 2 \Psi_* \right) d + \tfrac{\sigma}{\kmax} \Psi_* d \\ 
    \tfrac{- \sigma }{\kmax} \Psi_* d &> \tfrac{r}{\kmax} \left( \tfrac{V}{\mk} - 2 \Psi_* \right) d - \alpha , \\
    \sigma &< \tfrac{- r}{\Psi_*} \left( \tfrac{V}{\mk} - 2 \Psi_* \right) + \tfrac{\alpha \kmax}{\Psi_* d} = 2 r - \tfrac{r V}{\mk \Psi_*} + \tfrac{\alpha \kmax}{\Psi_* d} 
\end{align*}
and the proof is completed. 
\end{pf}
It should be stressed that this result is in agreement with Theorem 3, from \cite{dario}. Interestingly enough, our result has been found by starting form Greshgorin's circle theorem, while in \cite{dario} they analyse a $2 \times 2$ matrix by considering its trace and determinant. We have thus derived a similar result using different methodology.

\section{Symmetric Equilibrium in Structured Environment}\label{sec:regstruc}
In this section, we will assume that the equilibrium we reached is a symmetric equilibrium, cf. case 1 in Theorem 2. In such case, we have that $x_i^k = y_i^k = \xi_i^k = \mu_i^k$. Since this is a consensus equilibrium, by which we mean that the value is equal in every population $i$, this value is in fact equal to $\xi^k$, where the equilibrium value can still depend on the node degree $k$. We remark that for these symmetric consensus equilibria, $\theta_i^x = \theta_i^y = \theta$. Hence, 
\[ \sum_{j \in \mathcal{N}_i} \theta_j^x = \sum_{j \in \mathcal{N}_i} \theta_j^y = d \theta , \] 
where $d$ denotes the degree of the regular network. The dynamics for $x_i^k$ and $y_i^k$ are now described by 
\begin{align*}
    \dot{x}_i^k &= ( 1 - x_i^k - y_i^k ) \left( r \psi_k d \theta  + \gamma  \right) - x_i^k \left( \alpha + \sigma \psi_k d \theta  \right) , \\ 
    \dot{y}_i^k &= ( 1 - x_i^k - y_i^k ) \left( r \psi_k d \theta + \gamma  \right) - y_i^k \left( \alpha + \sigma \psi_k d \theta  \right) .
\end{align*}    
Equivalently we can write these dynamics in matrix form to get \footnotesize 
\begin{align}\begin{split} \label{eq:regularstrucdyn} 
    \begin{pmatrix} \dot{x}_i^k \\ \dot{y}_i^k \end{pmatrix} &= \begin{bmatrix} - (r + \sigma) \psi_k d \theta - \alpha - \gamma & - r \psi_k d \theta - \gamma \\ - r \psi_k d \theta - \gamma & - (r + \sigma) \psi_k d \theta - \alpha - \gamma \end{bmatrix} \begin{pmatrix} x_i^k \\ y_i^k \end{pmatrix} 
    \\ &\quad + \begin{pmatrix} r \psi_k d \theta + \gamma \\ r \psi_k d \theta + \gamma \end{pmatrix} .  \end{split} 
\end{align}
\normalsize The remainder of this section focuses on the study of the above system. By direct computation, one can show that the above affine system is stable by computing the eigenvalues of the system matrix.          
\begin{thm} The above system has eigenvalues $\lambda =  - \sigma \psi_k d \theta - \alpha$ and $\lambda =  - (2r + \sigma ) \psi_k d \theta - 2 \gamma - \alpha $, which are strictly negative and hence the system is asymptotically stable. Furthermore convergence is faster as $\psi_k$ increases.  \end{thm} 
Next, we study the equilibrium of system \eqref{eq:regularstrucdyn}. 
\begin{thm} For the system described by \eqref{eq:regularstrucdyn}, a symmetric equilibrium is given by \begin{align*} x_i^{k*} = y_i^{k*} = x^* = \dfrac{r \psi_k d \theta + \gamma }{(2 r + \sigma ) \psi_k d \theta + 2 \gamma + \alpha } . \end{align*} \end{thm} 
\begin{pf}  As the dynamics of \eqref{eq:regularstrucdyn} are of the form $\dot{x} = A x + c$, the equilibrium is simply given by $x^* = - A^{-1} c$. We remark that the determinant of $A$ is given by 
\[ \text{det}(A) = \Delta = 2 (\rpdt + \ga )(\spdt + \al ) + (\spdt + \al)^2 ,\]
which is nonzero as the system parameters are all positive. So matrix $A$ is always nonsingular. Performing straightforward computations to obtain $-A^{-1} c$, yield the desired result.
\end{pf}
One can show that if $\tfrac{\al}{\ga}> \tfrac{\sigma}{r}$, we have that $\tfrac{\partial x_i^{k*}}{\partial \psi_k} (\psi_k) > 0$ for all $\psi_k$, thus the equilibrium value $x_i^{k*}$ is increasing as the connectivity $\psi_k$ is increasing. In such case, the fraction of uncommitted individuals (within a cluster) is decreasing.

\section{Numerical Studies}\label{sec:simul}
\begin{figure}[ht!]
    \centering
    \includegraphics[scale=0.50]{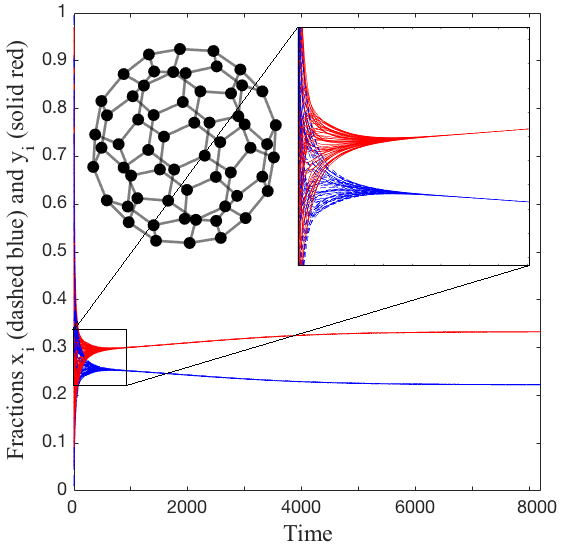}
    \caption{Dynamics on the regular networked Buckminster dome.}
    \label{fig:buckysimulg}
\end{figure}

Numerous simulations are run to support our findings. First, a simulation is performed to show the validity of Theorem \ref{thm:conseq}. We look at a regular and unweighted network given by the Buckminster Fuller geodesic dome, which is a regular graph on $60$ nodes where each node has degree~$3$, see the top left corner of Figure \ref{fig:buckysimulg}.  We performed the simulation and the results are shown in Figure \ref{fig:buckysimulg}. The values of the parameters were set as follows: $\gamma=0.2, \alpha=0.4, r=0.3$ and $\sigma=0.4$. 

The trajectories of $x_i$ are presented in dashed blue lines, while the trajectories of $y_i$ are presented in solid red lines. We observe that we reach a consensus equilibrium with $x_i = \xi = 0.222$ and $y_i = \mu =0.333$. We note that $\zeta = 1 - \xi - \mu = 0.444$. This is equal to $\frac{\alpha}{r d}$. Computing the values of $\xi$ and $\mu$ using the formulas in case 2 of Theorem \ref{thm:conseq} yield the same values as the results obtained from the simulation.

Next, we will perform simulations on the networked multi-population where the individual populations possess a structured environment. We assume that every node of the networked multi-population represents a complex network, and this complex network is clustered based on the connectivity.

As every population $i$ is now a complex network, cluster $k$ in population $i$ is simply the portion of nodes with $k$ connections. We have thus $k_{\text{max}}$ clusters, where the final cluster has the maximal connectivity $k_{\text{max}}$. We assume that every population of the networked multi-population has the same probability distribution of the complex network, and this discrete distribution follows a power law distribution. 

This complex network was constructed using the Barabasi-Albert model, making the complex network scale-free. We consider the networked multi-population where the overall network is given by the regular graph of the Buckminster dome. So, we have in fact an interconnected network of 60 complex networks, where each complex network is scale free and follows the degree distribution of the Barabasi-Albert model. 

\begin{figure}[ht!]
    \centering
    \includegraphics[scale=0.33]{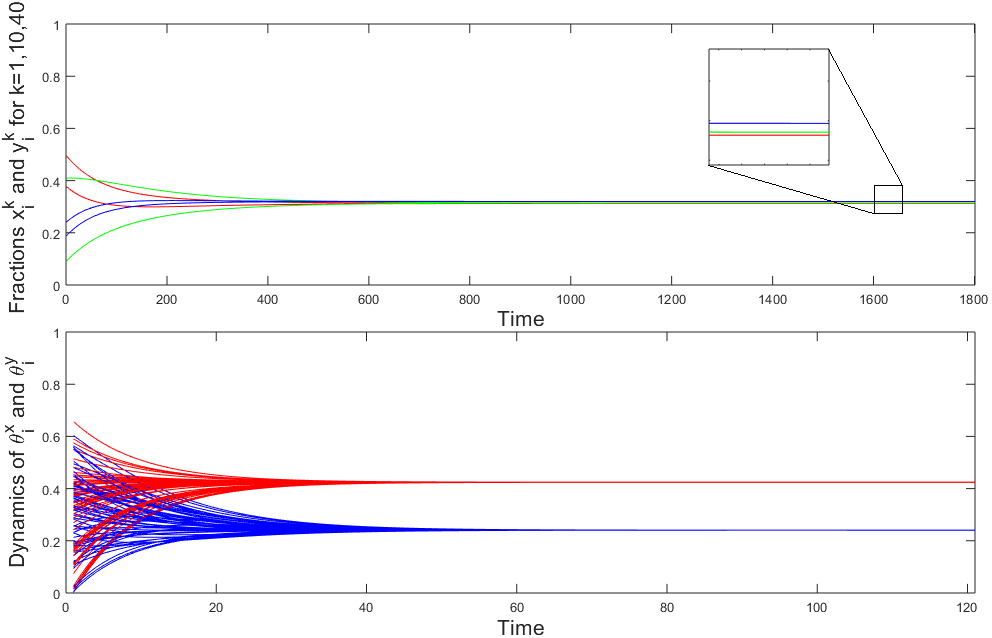}
    \caption{Top: evolutionary dynamics for different clusters of a population, with connectivity $1$, $10$ and $40$ respectively. Bottom: dynamics of all $\tix$ and $\tiy$ for every $i=1,2,\ldots,60$. }
    \label{fig:clusterdynam}
\end{figure}

We set the system parameters equal to $\gamma = 0.5, \alpha = 0.6, r=0.4$ and $\sigma = 0.3$. Then, we simulate the normalized dynamics of each cluster according to equation \eqref{eq:micromod}. At every iteration, for every population, we compute the $\theta_i^x$ and $\theta_i^y$. We have plotted the dynamics for different clusters $k=1,10,40$ of an arbitrary population of the multi-population. The results are given in Figure \ref{fig:clusterdynam} (top). Here, the normalised fractions $x_i^1$ and $y_i^1$ are plotted in red, $x_i^{10}$ and $y_i^{10}$ are plotted in green, and $x_i^{40}$ and $y_i^{40}$ are plotted in blue. We note that the six trajectories converge to three values: $x_i^1$ and $y_i^1$ converge to $0.3127$, $x_i^{10}$ and $y_i^{10}$ converge to $0.3143$ and $x_i^{40}$ and $y_i^{40}$ converge to $0.3190$. We remark that for every cluster, we reach a symmetric equilibrium since $x_i^k$ and $y_i^k$ converge to the same value. We note that this equilibrium value increases as the connectivity increases, which is in accordance with Theorem 6.

Finally, we perform a simulation on the model described by~\eqref{eq:micmac}. As starting point we consider again the multi-population of the Buckminster dome, and the values of $\sigma, r, \alpha$ and $\gamma$ are the same as of the previous simulation. As starting point of our simulation, we take the steady state values of $x_i^{k*}$ and $y_i^{k*}$ from the previous simulation, so that we can compute the second moments $\pix, \piy$ and $\piz$. Having fixed these, we now consider the dynamics given by \eqref{eq:micmac}. The time evolution of $\theta_i^x$ and $\theta_i^y$ are given in Figure \ref{fig:clusterdynam} (bottom), where $\tix$ is plotted in red and $\tiy$ is plotted in blue.

We observe that the $\tix$ and $\tiy$ reach a consensus equilibrium. It must be stressed that for these values of the system, the condition \eqref{eq:stabcond} is satisfied and we have reached asymptotic stability, validating Corollary 5.

\section{Conclusion}\label{sec:conclusion}
In this paper we studied bio-inspired evolutionary dynamics on a regular network representing a multi-population. The dynamics of an individual of each group now not only changes due to changes within its own group, but it takes into account the states of neighboring populations as well. First we gave a description of the average behavior of each population, and the states are the fractions of the population committed to each option. We showed that the steady state that is reached is a consensus equilibrium, and a sufficient condition for stability of this equilibrium in terms of a lower bound on the cross-inhibitory signals was given. Secondly, we added a structured environment to each population by assuming that each group was represented by a complex network. Within a group, it is then possible to cluster the group based on the internal connectivity of a player, and we derived a model for the dynamics of each cluster of each population. By aggregating the equations over all clusters of a single population, we found a description in terms of second moments. We analyzed this model and we presented stability results of the equilibrium. The paper is concluded by simulations validating our theoretical results. 

For future research, we will study a general network of a multi-population, where the condition that the underlying network is regular is dropped. We will also investigate the heterogeneous situation, where the system parameters are different for each group. Finally, the case of asymmetric parameters is of interest as well.


\bibliography{mybibfile}

\begin{thebibliography}{10}
\expandafter\ifx\csname url\endcsname\relax
  \def\url#1{\texttt{#1}}\fi
\expandafter\ifx\csname urlprefix\endcsname\relax\def\urlprefix{URL }\fi
\expandafter\ifx\csname href\endcsname\relax
  \def\href#1#2{#2} \def\path#1{#1}\fi

\bibitem{dario}
L.~Stella, D.~Bauso, Bio-inspired evolutionary dynamics on complex networks
  under uncertain cross-inhibitory signals, Automatica 100 (2019) 61--66.

\bibitem{j13}
A.~Pluchino, V.~Latora, A.~Rapisarda, Compromise and synchronization in opinion
  dynamics, The European Physical Journal B-Condensed Matter and Complex
  Systems 50~(1-2) (2006) 169--176.

\bibitem{hill}
S.~Tan, J.~Lu, G.~Chen, D.~J. Hill, When structure meets function in
  evolutionary dynamics on complex networks, IEEE Circuits and Systems magazine
  14~(4) (2014) 36--50.

\bibitem{bauso3}
L.~Stella, D.~Bauso, Evolutionary game dynamics for collective decision making
  in structured and unstructured environments, IFAC-PapersOnLine 50~(1) (2017)
  11914--11919.

\bibitem{bee1}
N.~F. Britton, N.~R. Franks, S.~C. Pratt, T.~D. Seeley, Deciding on a new home:
  how do honeybees agree?, Proceedings of the Royal Society of London. Series
  B: Biological Sciences 269~(1498) (2002) 1383--1388.

\bibitem{bee2}
R.~Gray, A.~Franci, V.~Srivastava, N.~E. Leonard, Multiagent decision-making
  dynamics inspired by honeybees, IEEE Transactions on Control of Network
  Systems 5~(2) (2018) 793--806.

\bibitem{bee3}
V.~Srivastava, N.~E. Leonard, Bio-inspired decision-making and control: From
  honeybees and neurons to network design, in: 2017 American Control
  Conference, IEEE, 2017, pp. 2026--2039.

\bibitem{multipoppaper}
D.~Bauso, Consensus via multi-population robust mean-field games, Systems \&
  Control Letters 107 (2017) 76--83.

\bibitem{opiniondyn}
R.~Hegselmann, U.~Krause, Opinion dynamics and bounded confidence models,
  analysis, and simulation, Journal of artificial societies and social
  simulation 5~(3).

\bibitem{ming}
W.~Yu, G.~Chen, M.~Cao, Consensus in directed networks of agents with nonlinear
  dynamics, IEEE Transactions on Automatic Control 56~(6) (2011) 1436--1441.

\bibitem{bullogoed}
W.~Mei, S.~Mohagheghi, S.~Zampieri, F.~Bullo, On the dynamics of deterministic
  epidemic propagation over networks, Annual Reviews in Control 44 (2017)
  116--128.

\bibitem{basar}
J.~Liu, P.~E. Par{\'e}, A.~Nedi{\'c}, C.~Y. Tang, C.~L. Beck, T.~Ba{\c{s}}ar,
  Analysis and control of a continuous-time bi-virus model, IEEE Transactions
  on Automatic Control 64~(12) (2019) 4891--4906.

\bibitem{onzeconf}
W.~Baar, D.~Bauso, Networked bio-inspired evolutionary dynamics on a
  multi-population, in: 18th European Control Conference, IEEE, 2019, pp.
  1023--1028.

\bibitem{frabl}
F.~Blanchini, S.~Miani, Set-theoretic methods in control, Springer, 2008.

\bibitem{thepaperonmydesk}
Y.~Moreno, R.~Pastor-Satorras, A.~Vespignani, Epidemic outbreaks in complex
  heterogeneous networks, The European Physical Journal B-Condensed Matter and
  Complex Systems 26~(4) (2002) 521--529.

\bibitem{thetaref}
Y.-Y. Ahn, H.~Jeong, N.~Masuda, J.~D. Noh, Epidemic dynamics of two species of
  interacting particles on scale-free networks, Physical Review E 74~(6) (2006)
  066113.

\end{thebibliography}

\end{document}